\documentclass[11pt]{amsart}

\newtheorem{theorem}{Theorem}[section]
\newtheorem{lemma}[theorem]{Lemma}

\theoremstyle{definition}
\newtheorem{definition}[theorem]{Definition}
\newtheorem{example}[theorem]{Example}

\theoremstyle{remark}
\newtheorem{remark}[theorem]{Remark}

\def\B{{\mathcal B}}

\def\P{{\bf P}}
\def\H{{\bf H}}
\def\R{{\bf R}}
\def\E{{\bf E}}
\def\I{{\bf I}}
\def\eps{{\varepsilon}}
\def\complexity{{\hbox{\rm complex}}}
\def \endprf{\hfill  {\vrule height6pt width6pt depth0pt}\medskip}
\def\emph#1{{\it #1}}
\def\textbf#1{{\bf #1}}

\begin{document}

\title{Szemer\'edi's regularity lemma revisited}

\author{Terence Tao}
\address{Department of Mathematics, UCLA, Los Angeles CA 90095-1555}
\email{tao@@math.ucla.edu}
\thanks{The author thanks Fan Chung Graham for helpful comments, and Jozsef Solymosi for encouraging the creation of this manuscript.  The author is also indebted to the anonymous referees for many useful suggestions and corrections.  The author is supported by a grant from the Packard Foundation.}

\vspace{-0.3in}
\begin{abstract}
Szemer\'edi's regularity lemma is a basic tool in graph theory, and also plays an important
role in additive combinatorics, most notably in proving Szemer\'edi's theorem on arithmetic 
progressions \cite{szemeredi}, \cite{szemeredi-4}.  In this note we revisit this lemma from the perspective of
probability theory and information theory instead of graph theory, and observe a slightly stronger variant of this 
lemma, related to similar strengthenings of that lemma in \cite{af}.  This stronger version of the regularity lemma was extended in \cite{tao:hyper} to reprove
the analogous regularity lemma for hypergraphs.
\end{abstract}

\maketitle

\section{Introduction}

Szemer\'edi's regularity lemma, introduced by Szemer\'edi in \cite{szemeredi}, is a fundamental tool in graph theory, and more precisely
in the theory of very large, dense graphs.  Roughly speaking, it asserts that given any such large dense graph $G$, and given an error tolerance
$0 < \eps \ll 1$, one can approximate $G$ by a much simpler object, namely a partition of the vertex set into $O_\eps(1)$ classes, together
with some edge densities between atoms of this partition, such that the approximation is ``$\eps$-regular'' on most pairs of this partition;
we will formalize these notations shortly.  This lemma can thus be viewed as a structure theorem for large dense graphs, approximating such
graphs to any specified accuracy by objects whose complexity is bounded independently of the number of vertices in the original graph.

The regularity lemma has had many applications in graph theory, computer science, discrete geometry and in additive combinatorics, see \cite{komlos} for a survey.  In particular, this lemma
and its variants play an important role in Szemer\'edi's celebrated theorem \cite{szemeredi} that any subset of the integers of positive
density contain arbitrarily long arithmetic progressions.  A variant of this structure theorem (also borrowing heavily from ideas in ergodic theory)
was also crucial in showing in \cite{gt-prime} that the primes contained arbitrarily long arithmetic progressions.  The lemma 
has also had a number of generalizations to hypergraphs of varying degrees of strength, see \cite{chung-hyper},
\cite{frankl}, \cite{frankl02}, \cite{nrs}, \cite{rs}, \cite{rodl}, \cite{gowers-hyper}, \cite{tao:hyper}.  The more recent formulations of the hypergraph lemma are in fact strong enough to rather easily imply Szemer\'edi's theorem on arithmetic progressions, as well as a multidimensional version due to Furstenberg and
Katznelson \cite{fk}.  They were also used in the recent paper \cite{tao-multiprime} establishing infinitely many constellations of any given shape in the Gaussian primes.

The proof of Szemer\'edi's lemma is now standard in the literature.  However, this standard proof is difficult to extend to the hypergraph
case; a direct application of the argument does give fairly easily a regularity lemma for hypergraphs (see  
\cite{chung-hyper}, \cite{frankl}), but that lemma does not seem to be strong enough for applications such as 
Szemer\'edi's theorem or the  Furstenberg-Katznelson theorem\footnote{The difficulty is that in the hypergraph situation, there are several levels of regularity or discrepancy that need to be controlled in order to yield a useful bound for arithmetic progressions or similar structures, and the lemma in \cite{chung-hyper} or \cite{frankl} controls only one of these discrepancies.  Later regularity lemmas control all of the relevant discrepancies, but there are some non-trivial technical issues concerning the relative sizes of the error estimates, as certain losses coming from one level of approximation must be compensated for by gains from the discrepancy bounds in other levels of approximation.}, except when concerning progressions or constellations consisting of at most three points (see \cite{soly-roth}).

In this paper we shall present a slightly different way of looking at Szemer\'edi's regularity lemma, which we used in \cite{tao:hyper} to obtain a hypergraph regularity lemma with sufficient strength for applications to Szemer\'edi-type theorems.  In this new perspective, one views the regularity lemma not as a structure theorem for large dense graphs, but rather
as a structure theorem for events or random variables in a product probability space.  This change of 
perspective is analogous to Furstenberg's highly successful approach to Szemer\'edi's theorem in \cite{furst}, in which the purely combinatorial result of Szemer\'edi was recast as
a statement about recurrence for arbitrary events or random variables in a probability-preserving system.  Just as Furstenberg's
change of perspective allowed the powerful techniques of ergodic theory to be brought to bear on the problem, the change of perspective here allows one to employ tools from probability theory and information theory to clarify the regularity lemma.  In particular we will use three very useful
concepts from those theories, namely \emph{$\sigma$-algebras} (partitions), \emph{conditional expectation} (relative density), and
\emph{entropy} (complexity).  As the parenthetical comments suggest, each of these concepts has a combinatorial analogue, however the author believes that there is some conceptual advantage to be gained by using a probabilistic and information-theoretic perspective rather than a graph-theoretic
one\footnote{The situation is somewhat analogous to that of the probabilistic method in combinatorics.  While every probabilistic argument could, in principle, be written in a deterministic way (replacing expectations by averages, etc.), it is undeniable that there are significant conceptual benefits in using a ``probabilistic way of thinking'' to approach combinatorial problems.}.  One byproduct of this new perspective is that
one discovers a stronger and more flexible version of the regularity lemma hiding underneath the standard one.  This stronger
version is difficult to state here without the requisite notational setup, but let us just say for now that it is closely related to a similar
improvement of the regularity lemma discovered recently\footnote{Note added in proof: a closely related version of this lemma was recently introduced in \cite{af}, \cite{alonshapira}.  See also \cite{lovasz-sz} for yet another perspective on the regularity lemma, this time from functional analysis.}  in \cite{rs}, in which it was observed that the regularity of the large dense graph $G$
relative to the partition given by that lemma can be vastly improved after adding or removing a small number of edges from $G$.  This strengthened
version of the regularity lemma turns out to be quite amenable for iterating, and thus gives a relatively painless proof of the hypergraph regularity lemma; see \cite{tao:hyper}.

We will turn to the details in later sections, but for now let us just give an informal discussion which already shows that the regularity lemma can be
viewed in information theoretic terms rather than graph theoretic terms.  It will be convenient to work with bipartite graphs.
Let $G = (V_1,V_2,E)$ be a large dense bipartite graph.  Let $x_1$ and $x_2$ be two vertices
selected independently and uniformly at random from $V_1$ and $V_2$ respectively; thus $x_1$ and $x_2$ are independent random variables, taking values in $V_1$ and $V_2$ respectively.  The edge set
$E$ can now be re-interpreted as a probabilistic event, namely the event that the pair $(x_1,x_2)$ lies in $E$.  We shall abuse notation and refer to
this event also as $E$, thus $E$ is now some event determined by the random variables $x_1,x_2$ (or more precisely, it lies in the $\sigma$-algebra
generated by the random variables $x_1$ and $x_2$).  Many of the important statistics about the edge set $E$ can now be recast in terms of the
event $E$; for instance, the edge density of the edge set $E$ is equal to the probability of the event $E$, or equivalently the expectation of
the indicator random variable $1_E$.  Similarly one can view relative edge densities of $E$ as conditional expectations of $1_E$.  

We have already observed that $E$ is, in principle, determined by $x_1$ and $x_2$.  However, from an information-theoretic perspective
this determinism relationship can be very ``high-complexity'' or ``fine-scaled'', in a sense we shall describe shortly.  
If the vertex sets $V_1,V_2$ have $N$ elements, then the random variables $x_1$ and $x_2$ have
a Shannon entropy of $\log_2 N$ (they can be described by roughly $\log_2 N$ bits each). On the other hand, the event $E$ (or the Boolean function $1_E$) has a Shannon entropy of at most $\log_2 2 = 1$ (it can be described by one bit).  If $N$ is very
large, we thus see that there is much more information contained in the random variables $x_1$ and $x_2$ than is contained in the event $E$.  To put it
another way, knowing that the event $E$ is true or false (i.e., that the pair $(x_1,x_2)$ is an edge in $G$ or not) does not even begin to let one
determine the exact values of $x_1$ and $x_2$.  Indeed, in the extreme case when the graph $G$ is a random (or pseudorandom) graph, the event $E$
behaves almost as if it were \emph{independent} of the random variables $x_1$ and $x_2$, despite being actually determined by these variables.  More precisely, if $A_1$ is any event determined by $x_1$ (thus $A_1$ can be thought of as the event that $x_1$ lies in a fixed subset of $V_1$, which by abuse of notation we shall also call $A_1$), and $A_2$ is any event determined by $x_2$, then in the random or pseudorandom case the event $E$ will be
almost completely uncorrelated with the events $A_1,A_2$.  This corresponds to the well-known fact that when $G$ is a random or pseudorandom graphs,
the relative edge density between two large sets $A_1,A_2$ in $V_1, V_2$ will, with high probability, be very close to the global edge density of $G$.
(Note that if $A_1$ and $A_2$ were small sets, i.e. events of very low probability, then the correlation, or more precisely the \emph{mutual information}, with $E$ would automatically be small.)

Let us summarize the above discussion in information-theoretic terms.  If one is given all $\log_2 N$ bits of $x_1$, and all $\log_2 N$ bits of $x_2$,
then the single-bit event $E$ is completely determined.  But if $G$ is random or pseudorandom, and one is only given one bit 
of $x_1$ (specifically, whether $x_1$ lies in a fixed set $A_1$) and one bit of $x_2$, one learns almost no information about the bit $E$.  Let us informally describe this by saying that $E$ is approximately independent of $x_1$ and $x_2$ at ``coarse scales'' - when only a few bits of $x_1$ and $x_2$ are known, even though $E$ is determined by $x_1$ and $x_2$ at ``fine scales'' - when most or all of the bits of $x_1$ and $x_2$ are known.

Of course, if $G$ is not pseudorandom, then $E$ can be highly correlated with a few special bits of $x_1$ and $x_2$.  To take an extreme opposite case to the pseudorandom case, suppose that $G$ is a complete bipartite graph connecting all the vertices of a set $A_1 \subseteq V_1$ to that of a 
set $A_2 \subseteq V_2$, and not connecting any other pairs of vertices. Then the event $E$ is completely determined by one bit of $x_1$ (namely, whether it lies in $A_1$) and one bit 
of $x_2$ (namely, whether it lies in $A_2$).

Furthermore, it is possible for $G$ to be a hybrid between these two extremes.  Suppose now that $G$ is a pseudorandom subgraph of the complete bipartite graph connecting $A_1$ to $A_2$.  Then $E$ is no longer determined by the one special bit of $x_1$ associated to $A_1$, and the one special bit of $x_2$ associated to $A_2$.  However, it is now approximately independent at coarse scales of $x_1$ and $x_2$, \emph{conditioning on $A_1$ and $A_2$}.  In other words, once the events $A_1$ and $A_2$ are known to be true or false, the event $E$ is then approximately independent to any further bits of information arising from $x_1$ and $x_2$.  In graph theory terms, this means that when restricting $V_1$ to $A_1$ or its complement,
and restricting $V_2$ to $A_2$ or its complement, the induced subgraph of $G$ behaves pseudorandomly (with some edge density depending on which
sets were being restricted to).

The information-theoretic version of the Szemer\'edi regularity lemma is an assertion, roughly speaking, that every event $E$ is a hybrid
of the two extremes in the sense given above.  Very informally, given
any two high-entropy random variables $x_1$ and $x_2$, and given any event $E$, it is possible to find some low-entropy random variable $Z_1$ determined by $x_1$, and a low-entropy random variable $Z_2$ determined by $x_2$, such that $E$ is approximately independent of $x_1$ and $x_2$
conditioning on $Z_1$ and $Z_2$.  Again being very informal, this means that there exist a small number of bits from $x_1$ and $x_2$
which correlate with $E$, and such that no further bits from $x_1$ and $x_2$ have much of a correlation with $E$.
Interestingly, this formulation of the regularity lemma requires no independence properties of $x_1$ and $x_2$, and also does not require $E$
to be determined by $x_1$ and $x_2$; but we do not know any applications of this more general version.  

One can view the low-entropy random variables $Z_1$, $Z_2$ discussed above as ``approximations'' to the event $E$, where the approximation is in some
coarse information-theoretic sense.  It turns out that the proof of the regularity lemma (see Lemma \ref{mainlemma} below) in fact yields \emph{two} such approximations, a
``coarse approximation'' $Z_1, Z_2$ and a ``fine approximation'' $Z'_1, Z'_2$.  The coarse approximation has low entropy.  The fine approximation
has significantly higher entropy, but it is an exceedingly accurate approximation to $E$; in particular, any error arising from this approximation
can exceed any losses coming from the entropy of the coarse approximation, in a way which can be made precise using a ``growth function'' $F: \R^+ \to \R^+$.  Finally, the coarse and fine approximations will be close to each other, both in an $L^2$ sense, and also in an information theoretic 
sense.  We will make these statements more precise later, however we remark for now that the presence of the new parameter $F$, used to compare
the accuracy of the fine approximation against the entropy of the coarse approximation, is very suitable for iteration purposes, and allows
one to extend the regularity lemma to the hypergraph setting, in which one has multiple random variables $x_1,\ldots,x_d$ instead of just two, and furthermore one is interested in low-entropy approximations to an event which arise not only from individual random variables $x_i$,
but also from joint random variables such as $(x_i, x_j)$ (and the approximations coming from the joint random variables should themselves be approximated by other, lower-order random variables).  See \cite{tao:hyper}.  A closely related regularity lemma,
which also involves an arbitrary growth function $F$, has also recently appeared in \cite{af} in applications to property testing.

\section{A probabilistic formulation}

Before we give the rigourous information-theoretic version of the Szemer\'edi regularity lemma, let us first give a standard formulation of the
lemma, and also a probabilistic formulation which can be viewed as an intermediate formulation bridging the graph-theoretic version and
the information-theoretic\footnote{We say a formulation is ``probabilistic'' if it involves such concepts as probability spaces, $\sigma$-algebras, random variables, (conditional) expectation, and correlation.  We say a formulation is ``information-theoretic'' if it involves such concepts as probability spaces, $\sigma$-algebras, random variables, (conditional) entropy, and mutual information.  Clearly these two perspectives share much in common, for instance the concept of independence is important in both.} version of the lemma.  We begin with the graph-theory version; again, it is convenient to restrict ones attention to bipartite graphs.

We use $O(X)$ to denote any quantity bounded in magnitude by $CX$ for some absolute constant $C > 0$, and more generally
we use $O_{a_1,\ldots,a_k}(X)$ to denote any quantity bounded in magnitude by $C(a_1,\ldots,a_k) X$, where $C(a_1,\ldots,a_k) > 0$
depends on the parameters $a_1,\ldots,a_k$. We also use $|A|$ to denote the cardinality of a finite set $A$.

\begin{definition}  A \emph{bipartite graph} is a triplet $(V_1,V_2,E)$ where $V_1,V_2$ are two finite non-empty sets,
and $E \subset V_1 \times V_2$.  If $\eps > 0$, we say that
a bipartite graph $(V_1,V_2)$ is \emph{$\eps$-regular} if we have
\begin{equation}\label{ecap}
|E \cap (A_1 \times A_2)| = \frac{|A_1 \times A_2|}{|V_1 \times V_2|} |E| + O( \eps |V_1 \times V_2| )
\end{equation}
for all $A_1 \subseteq V_1$ and $A_2 \subseteq V_2$.  
\end{definition}

\begin{remark} While we assert that \eqref{ecap} holds for all subsets $A_1$, $A_2$ of $V_1, V_2$, this condition is only
non-trivial for large subsets; it holds trivially when $|A_1 \times A_2| = O( \eps |V_1 \times V_2| )$.  Thus this definition of $\eps$-regularity is essentially equivalent to other formulations of regularity in the literature in which a lower bound
is imposed on the size of $A_1$ and $A_2$.
\end{remark}

\begin{theorem}[Szemer\'edi regularity lemma, graph-theoretic version]\label{srl}  Let $(V_1,V_2,E)$ be a bipartite graph, and let $0 < \eps \leq 1$.
Assume that $V_1$ and $V_2$ are large depending on $\eps$, thus $|V_1|, |V_2| \geq O_\eps(1)$.
Then there exists a positive integer $J = O_\eps(1)$ and decompositions
$$ V_i = V_{i,0} \cup V_{i,1} \cup \ldots \cup V_{i,J}$$
for $i=1,2$ with the following properties:
\begin{itemize}
\item (Exceptional set) For all $i=1,2$, we have $|V_{i,0}| = O(\eps |V_i|)$.
\item (Uniform partition) For all $i=1,2$ and $1 \leq j \leq J$ we have $|V_{i,j}| = |V_{i,j'}|$.
\item (Regularity) The induced bipartite graph $(V_{1,j_1},V_{2,j_2}, E \cap (V_{1,j_1} \times V_{2,j_2}))$ is $\eps$-regular for all but $O(\eps J^2)$ of the pairs $1 \leq j_1 \leq M$, $1 \leq j_2 \leq J$.
\end{itemize}
\end{theorem}

\begin{remark} The bound $J = O_\eps(1)$ is a little deceptive, as it conceals the fact that $J$ can in fact be extremely large depending on 
$1/\eps$, indeed there are examples where $J$ grows like an exponential tower of height equal to some power of $1/\eps$ (see \cite{gowers-sz}).
However, the key point is that the bound on $J$ does not depend on the cardinality of $V_1$ or $V_2$.  Indeed we shall shortly give a probabilistic formulation in which $V_1$ and $V_2$ could be infinite (cf. \cite{lovasz-sz}).
\end{remark}

We now give a probabilistic generalization of the above regularity lemma.  We first recall some standard notation from probability theory.

\begin{definition}[Probability space]  
A \emph{probability space} is a triple $(\Omega, \B_{\operatorname{max}}, \P)$, where $\Omega$ is a set (called the \emph{sample space}),
$\B_{\operatorname{max}}$ is a $\sigma$-algebra\footnote{A $\sigma$-algebra is a collection $\B$ of sets in the probability space $\Omega$ which is closed under (countable) unions, intersections, and complements, and contains the empty set and $\Omega$.  In the our applications $\B$ will typically be finite, in which case it can be identified with a finite partition $\Omega = \Omega_1 \cup \ldots \cup \Omega_M$ of the underlying probability space.  Indeed, the cells of this partition are the atoms (minimal non-empty elements) of $\B$, while $\B$ itself consists of all the sets which are unions of zero or more atoms in the partition.} of sets of $\Omega$ (the elements of $\B_{\operatorname{max}}$ being the \emph{events}), and $\P$ is a probability measure on 
$\B_{\operatorname{max}}$ (thus it is non-negative and has total mass one).  A \emph{random variable} is any measurable map $X: \Omega \to K$ to some space $K$
(which will typically either be a finite set, or the real line).   We let $L^1(\B_{\operatorname{max}})$ denote the space of real-valued, absolutely integrable random variables; as is customary we identify two random variables if they agree outside of an event of zero probability.  
If $X \in L^1(\B_{\operatorname{max}})$, we let $\E(X)$ denote the expectation of $X$.  In particular, if $E$ is an event, then $\E(1_E) = \P(E)$.
\end{definition}

\begin{remark}
For application to the regularity lemma, $\Omega$ will be a finite set, and $\B_{\operatorname{max}}$ will be the algebra of all subsets of $\Omega$, so there will
be no issues as to whether a random variable is measurable or integrable.  However, it is interesting to note that the arguments we give below extend
with no difficulty whatsoever to the case of infinite probability spaces.
\end{remark}

\begin{example}\label{running-1} Our primary application will be to bipartite graphs, say between two vertex classes $V_1$ and $V_2$.  In this case we can take $\Omega = V_1 \times V_2$, $\B_{\operatorname{max}}$ to be the power set of $\Omega$ (thus all subsets of $\Omega$ are measurable events), and $\P$ to be the uniform probability measure on $\Omega$; this corresponds to the operation of sampling two vertices $x_1$ and $x_2$ uniformly and independently at random from $V_1$ and $V_2$ respectively.  In this case, all functions $X: V_1 \times V_2 \to \R$ are measurable, and the expectation is just the average value on $V_1 \times V_2$.
\end{example}

A crucial concept from probability theory is that of \emph{conditional expectation}.

\begin{definition}[Conditional expectation]  Let $(\Omega, \B_{\operatorname{max}}, \P)$ be a probability space, and let $\B$ be a sub-$\sigma$-algebra of $\B_{\operatorname{max}}$.
If we let $L^2(\B)$ be the Hilbert space of $\B$-measurable, square-integrable real-valued random variables, with the usual norm
$\| X \|_{L^2(\B)} := \E( |X|^2 )^{1/2}$, then $L^2(\B)$ is a closed subspace of
$L^2(\B_{\operatorname{max}})$, and we let $X \mapsto \E(X|\B)$ be the associated orthogonal projection map from $L^2(\B_{\operatorname{max}})$ to $L^2(\B)$; thus for any square-integrable random variable $X \in L^2(\B_{\operatorname{max}})$, $\E(X|\B)$ will be a square-integrable $\B$-measurable random variable.
\end{definition}

The conditional expectation can be defined explicitly in the case when $\B$ is finite, which is in fact the only case we will need in this paper.
In such a case, the $\sigma$-algebra $\B$ is generated by a finite number of disjoint events $A_1, \ldots, A_n$ of positive probability, possibly together with some additional events of zero probability which we can safely ignore.  If $X \in L^2(\B_{\operatorname{max}})$, the conditional expectation $\E(X|\B)$
will be equal (almost surely) to $\E( X | A_i ) := \frac{1}{\P(A_i)} \E( X 1_{A_i} )$ on each event $A_i$.

Next, we define the complexity of a $\sigma$-algebra, which is a simplified version of the Shannon entropy.

\begin{definition}[Complexity]  Let $\B$ be a finite $\sigma$-algebra in a probability space $(\Omega, \B_{\operatorname{max}}, \P)$.  Then the
\emph{complexity} $\complexity(\B)$ of $\B$ is defined as the least number of events needed to generate $\B$ as a $\sigma$-algebra.
\end{definition}

Informally, a finite $\sigma$-algebra of complexity $M$ can be described using $M$ bits of information (equivalently, it contains at most $2^M$ atoms).

If $\B, \B'$ are two sub-$\sigma$-algebras of $\B_{\operatorname{max}}$, we let $\B \vee \B'$ denote the smallest $\sigma$-algebra which contains both
$\B$ and $\B'$.  Note that if $\B$ and $\B'$ are finite, then $\B \vee \B'$ is also finite, with the sub-additivity property
$$ \complexity(\B \vee \B') \leq \complexity(\B) + \complexity(\B').$$

\begin{example}\label{running-2} We continue the running example in Example \ref{running-1}. Any partition $V_1 = V_{1,1} \cup \ldots \cup V_{1,M}$ of the first vertex class induces a partition $V_1 \times V_2 = (V_{1,1} \times V_2) \cup \ldots \cup (V_{1,M} \times V_2)$ of the probability space $\Omega$ and hence creates a sub-$\sigma$-algebra $\B_1$ of $\B_{\operatorname{max}}$, which in information-theoretic terms captures the information of which cell of the partition the first vertex $x_1$ belongs to.  The complexity of $\B_1$ is essentially $\log_2 M$.  If we have another partition $V_2 = V_{2,1} \cup \ldots \cup V_{2,M}$ of the second vertex class we can form another $\sigma$-algebra $\B_2$, and thence create the joint $\sigma$-algebra $\B_1 \vee \B_2$, whose atoms
are pairs $V_{1,i} \times V_{2,j}$ and whose complexity is essentially $2\log_2 M$ (assuming for sake of discussion that all the cells in the partitions are non-empty).  
If $X: V_1 \times V_2 \to \R$ is any random variable (which one can think of as a weight
function assigning a number to each putative edge $(x_1,x_2)$), the conditional expectation $\E(X|\B_1 \vee \B_2)$ is then the function which on each pair of cells $V_{1,i} \times V_{2,j}$ takes a value equal to the relative density
$\frac{1}{|V_{1,i}| |V_{2,j}|} \sum_{x_1 \in V_{1,i}} \sum_{x_2 \in V_{2,j}} X(x_1,x_2)$ of $X$ on this pair of cells.  We remark
that when $X$ is the indicator function $X = 1_E$ of a graph, the $L^2$ norm of this conditional expectation (which we shall 
refer to here as the \emph{energy}) is a familiar concept in the standard treatment of the regularity lemma and is usually 
referred to as the \emph{index} of the partitions $\B_1, \B_2$.
\end{example}

We now give a probabilistic Szemer\'edi regularity lemma, which we state in considerably more generality than we need to establish
Theorem \ref{srl}.

\begin{theorem}[Szemer\'edi regularity lemma, probabilistic version]\label{measure-lemma}  Let $(\Omega, \B_{\operatorname{max}}, \P)$ be a probability space, let $(\B_{i,max})_{i \in I}$ be a finite collection of sub-$\sigma$-algebras of $\B_{\operatorname{max}}$, and let $X \in L^2(\B_{\operatorname{max}})$ be a random variable
with $\| X \|_{L^2(\B_{\operatorname{max}})} \leq 1$.
Let $\eps > 0$ be a number, let $m \geq 0$, and let $F: \R^+ \to \R^+$ be an arbitrary monotone increasing
function.  Then there exists finite sub-$\sigma$-algebras $\B_i \subseteq \B'_i \subseteq \B_{i,max}$ for each $i \in I$, and a non-negative real number\footnote{It may be helpful to the reader to think of $M$ as simply being the quantity
$\max_{i \in I}( m, \complexity(\B_i) )$.  Thus the upper bound on $M$ translates to an upper bound on the complexity of the coarse partitions $\B_i$, while the estimate \eqref{xbfine} asserts, roughly speaking, that the accuracy of the fine partitions exceeds the complexity of the coarse partitions (and also exceeds any specified constant $m$) by an arbitrary growth function $F$.}
 $M$, obeying the following bounds:
\begin{itemize}
\item (Size of $M$) We have $M \geq m$ and $M = O_{\eps, F, m}(1)$.
\item (Complexity bound) We have $\complexity(\B_i) \leq M$ for all $i \in I$.
\item (Coarse and fine approximations are close) We have
\begin{equation}\label{coarsefine}
 \left\| \E( X | \bigvee_{i \in I} \B'_i ) - \E( X | \bigvee_{i \in I} \B_i ) \right\|_{L^2(\B_{\operatorname{max}})} \leq \eps.
 \end{equation}
\item (Fine approximation is extremely accurate) For any collection $(A_i)_{i \in I}$ of events with $A_i \in \B_{i,max}$ for all $i \in I$,
we have
\begin{equation}\label{xbfine}
 \left|\E\left( \bigl(X - \E( X | \bigvee_{i \in I} \B'_i )\bigr) \prod_{i \in I} 1_{A_i} \right)\right| \leq \frac{1}{F(M)}.
\end{equation}
\end{itemize}
\end{theorem}

\begin{remark} In the application to Theorem \ref{srl}, we will only need this theorem in the special case when $X = 1_E$ is
an indicator function, when $I = \{1,2\}$, when $\B_{1,max}$, $\B_{2,\operatorname{max}}$ are finite and independent, with each atom having equal probability, 
$F$ is essentially the exponential function,
and $\B_{\operatorname{max}} = \B_{1,\operatorname{max}} \vee \B_{2,\operatorname{max}}$.  However the more general version above is no harder to prove than this special case.  One can
also generalize to the case when $X = (X_1,\ldots,X_n)$ is vector-valued, taking values in $\R^n$; on the graph level, this would correspond to
regularizing $n$ graphs simultaneously using a single partitioning of the vertex classes.  This vector-valued generalization is useful for iteration purposes, in order to easily obtain the corresponding hypergraph regularity lemma; this generalization is implicit in \cite{tao:hyper}.  
\end{remark}

\begin{remark}
Informally, this theorem starts with a square-integrable random variable $X$, 
and some reference $\sigma$-algebras $\B_{i,\operatorname{max}}$.  It then creates two approximations to $X$, namely a coarse approximation 
$\E(X| \bigvee_{i \in I} \B_i)$ and a fine approximation $\E(X| \bigvee_{i \in I} \B'_i)$.
The coarse approximation depends on only $M$ ``bits'' of information from each of the $\B_{i,\operatorname{max}}$, where $M$ is a quantity for which we have
some bounds.  The fine approximation is rather close to the coarse approximation in $L^2(\B_{\operatorname{max}})$ norm.  Finally, the fine approximation is
extremely accurate, in the sense that adding an additional bit of information from each of the $\B_{i,\operatorname{max}}$ can only create an additional correlation of at most $1/F(M)$, where $F(M)$ is a function of $M$ which can be specified in advance to be as rapidly growing as one pleases.
(Of course, there is a price to pay in selecting a function $F$ which grows too rapidly, which is that the upper bound on $M$ will deteriorate.)
Somewhat remarkably, no independence or dependence assumptions between $X$ and the $\B_{i,\operatorname{max}}$ need to be made in order for this theorem to be applicable.
\end{remark}

We will prove Theorem \ref{measure-lemma} in the next section.  For the remainder of this section, we show how Theorem \ref{measure-lemma} implies Theorem \ref{srl}.

\begin{proof}[Proof of Theorem \ref{srl} assuming Theorem \ref{measure-lemma}]  Let $G = (V_1,V_2,E)$ be a bipartite graph, thus $E$ can be
viewed as a subset of $V_1 \times V_2$.  We then define a probability space by setting the sample space $\Omega := V_1 \times V_2$,
setting the $\sigma$-algebra $\B_{\operatorname{max}} = 2^\Omega$ be the space of all subsets of $\Omega$, and setting $\P$ be the uniform probability
measure on $\Omega$.  In particular, $E$ is now an event in $\B_{\operatorname{max}}$.
As mentioned in the introduction, this probability space corresponds to the space generated by selecting vertices $x_1, x_2$
from $V_1, V_2$ independently and uniformly.  We then set $I := \{1,2\}$, and set $\B_{1,\operatorname{max}} := \{ A_1 \times V_2: A_1 \subseteq V_1 \}$
and $\B_{2,\operatorname{max}} := \{ V_1 \times A_2: A_2 \subseteq V_2 \}$, thus $\B_{1,\operatorname{max}}$ and $\B_{2,\operatorname{max}}$ are the $\sigma$-algebras 
generated by the random variables $x_1$ and $x_2$ respectively.  We set $X := 1_E$; clearly $\|X\|_{L^2(\B_{\operatorname{max}})} \leq 1$.

We now apply Theorem \ref{measure-lemma}, with the growth function $F: \R^+ \to \R^+$ to be chosen later, and $\eps$ replaced by
$\eps^{3/2}$.  This gives us some
$\sigma$-algebras $\B_1 \subseteq \B'_1 \subseteq \B_{1,\operatorname{max}}$ and $\B_2 \subseteq \B'_2 \subseteq \B_{2,\operatorname{max}}$ and a non-negative
quantity $M = O_{F,\eps}(1)$ such that
\begin{align}
\complexity(\B_1), \complexity(\B_2) &\leq M \label{complex-1} \\
\| \E( 1_E | \B'_1 \vee \B'_2 ) - \E( 1_E | \B_1 \vee \B_2 ) \|_{L^2(\B_{\operatorname{max}})} &\leq \eps^{3/2} \label{nearfine-1}\\
\left|\E\left( \bigl(1_E - \E( 1_E | \B'_1 \vee \B'_2 )\bigr) 1_{A_1 \times A_2} \right)\right| &\leq \frac{1}{F(M)}\label{finegood-1}
\end{align}
for all $A_1 \subseteq V_1, A_2 \subseteq V_2$.

Now let $J$ be a large integer to be chosen later; we will eventually show $J = O_\eps(1)$.  
By hypothesis we may take $|V_1|, |V_2| > J$.
For each $i \in \{1,2\}$, the finite $\sigma$-algebras $\B_i$ consists of at most $2^M$ 
atoms, thanks to \eqref{complex-1}.  Then we can subdivide each of these atoms arbitrarily into sets of size $\lfloor \frac{|V_i|}{(1+O(\eps)) J} \rfloor$, plus an error of size $O(|V_i|/J)$.  
Combining all of the errors into a single exceptional set $V_{i,0}$, we obtain a partition
$$ V_i = V_{i,0} \cup V_{i,1} \cup \ldots V_{i,J},$$
where the sets $V_{i,1},\ldots,V_{i,J}$ all have the same cardinality (comparable to $|V_i|/J$), 
and each lies in an atom of $\B_i$, and the exceptional set $V_{i,0}$ obeys the bounds
$$ |V_{i,0}| = O( \eps |V_i| ) + O( 2^M |V_i| / J ).$$
Thus, if we choose $J$ to be the nearest integer to $2^M / \eps$, we obtain $|V_{i,0}| = O( \eps |V_i| )$ as desired.  Also
we observe that since $M = O_{F,\eps}(1)$, we have $J = O_{F,\eps}(1)$.

Now consider an induced bipartite graph $G_{j_1,j_2} := (V_{1,j_1},V_{2,j_2}, E \cap (V_{1,j_1} \times V_{2,j_2}))$ where $1 \leq j_1,j_2 \leq J$.  
Suppose we wish to show that $G_{j_1,j_2}$ is $\eps$-regular, thus
$$
|E \cap (A_1 \times A_2)| = \frac{|E \cap (V_{1,j_1} \times V_{2,j_2})|}{|V_{1,j_1} \times V_{2,j_2}|} |A_1 \times A_2| + O( \eps |V_{1,j_1}| |V_{2,j_2}| )$$
whenever $A_1 \subseteq V_{1,j_1}$ and $A_2 \subseteq V_{2,j_2}$. By the triangle inequality (and by specializing the estimate below
to the case $A_1 = V_{1,j_1}, A_2 = V_{2,j_2}$), it suffices to find a quantity $d$ which is
independent of $A_1,A_2$ (but which depends on $E$, $V_{1,j_1}, V_{2,j_2}$) such that
$$
|E \cap (A_1 \times A_2)| = d |A_1 \times A_2| + O( \eps |V_{1,j_1}| |V_{2,j_2}| )$$
whenever $A_1 \subseteq V_{1,j_1}$ and $A_2 \subseteq V_{2,j_2}$.  Dividing by $|V_1| |V_2|$, we can rewrite this as
$$ \E( (1_E - d) 1_{A_1 \times A_2} ) = O( \eps / J^2 ).$$
Observe that $A_1 \times A_2$ is contained in a single atom of $\B_1 \vee \B_2$.  Thus we may take $d := \E( 1_E | \B_1 \vee \B_2 )$ on this
atom.  Our task is thus to establish
$$ \E( (1_E - \E( 1_E | \B_1 \vee \B_2)) 1_{A_1 \times A_2} ) = O( \eps / J^2 ).$$
From \eqref{finegood-1} we have
$$ \E( (1_E - \E( 1_E | \B'_1 \vee \B'_2 )) 1_{A_1 \times A_2} ) = O( 1 / F(M) )$$
and so if we choose $F(M) := 2^{2M} / \eps^3$ then we have
$$ \E( (1_E - \E( 1_E | \B'_1 \vee \B'_2 )) 1_{A_1 \times A_2} ) = O( \eps^3 / 2^{2M} ) = O( \eps / J^2 ).$$
Note that we now have $J = O_{F,\eps}(1) = O_\eps(1)$ as desired.  Thus, in order to establish $\eps$-regularity of
$G_{j_1,j_2}$, it suffices by the triangle inequality to establish that
$$ \E( |\E( 1_E | \B'_1 \vee \B'_2 ) - \E( 1_E | \B_1 \vee \B_2 )| 1_{V_{1,j_1} \times V_{2,j_2}} ) = O( \eps / J^2 ).$$
Note that $\E( 1_{V_{1,j_1} \times V_{2,j_2}} ) = O( 1 / J^2 )$.  Thus by Cauchy-Schwarz, it would thus suffice to show that
\begin{equation}\label{cz}
 \E( |\E( 1_E | \B'_1 \vee \B'_2 ) - \E( 1_E | \B_1 \vee \B_2 )|^2 1_{V_{1,j_1} \times V_{2,j_2}} ) = O( \eps^2 / J^2 ).
\end{equation}
On the other hand, from \eqref{nearfine-1} we have
$$ \E( |\E( 1_E | \B'_1 \vee \B'_2 ) - \E( 1_E | \B_1 \vee \B_2 )|^2 ) = O( \eps^3 ).$$
Thus there are at most $O(\eps J^2)$ pairs $(j_1,j_2)$ for which \eqref{cz} fails.  Thus we have
$\eps$-regularity for all but at most $O(\eps J^2)$ pairs, as desired.
\end{proof}

\begin{remark} It is clear from the argument that we can enforce a lower bound on the number $J$ of partitions, simply by
setting the parameter $m$ equal to a large number rather than equal to zero, since this will give a lower bound for $M$ and hence for $J$.
Of course, this will also increase the lower bound required for $|V_1|, |V_2|$, although in applications the cases when $|V_1|$ or $|V_2|$ are small
tend to be fairly easy (and the regularity lemma is of little use in such situations anyway).  Also, by considering multiple vertex sets $(V_i)_{i \in I}$ instead of just two, one can prove a version of hypergraph regularity lemma (similar to the early hypergraph
lemma in \cite{chung-hyper}) by a similar argument to the one given above; we omit
the details.  However to obtain the stronger and more modern versions of the hypergraph regularity lemma one needs to apply results such as the one above repeatedly; see \cite{tao:hyper} for more details.
\end{remark}

\section{Proof of Theorem \ref{measure-lemma}}

We now give the proof of Theorem \ref{measure-lemma}.  Let us fix $(\Omega, \B_{\operatorname{max}}, \P)$, $(\B_{i,\operatorname{max}})_{i \in I}$,
$X$, $\eps$, $m$, $F$.  A crucial concept in the proof (as in the standard proof of
the regularity lemma) will be that of the \emph{energy} (or index) of a $\sigma$-algebra (or partition).  This energy 
has a particularly simple description in the language of conditional expectation:

\begin{definition}  For any $\sigma$-algebra $\B \subseteq \B_{\operatorname{max}}$, we define the \emph{energy} ${\mathcal E}(\B)$ of $\B$ to be the
quantity
$$ {\mathcal E}(\B) := \| \E( X | \B ) \|_{L^2(\B_{\operatorname{max}})}^2.$$
Informally, ${\mathcal E}(\B)$ measures how close the subspace $L^2(\B)$ of the Hilbert space $L^2(\B_{\operatorname{max}})$ gets to containing the vector $X$.
\end{definition}

\begin{remark} In the running example of Example \ref{running-2}, with $X$ the indicator function of a graph and $\B = \B_1 \vee \B_2$, the energy corresponds to the \emph{index} of the partitions associated to $\B_1, \B_2$, as used for instance in \cite{szemeredi}.
\end{remark}

From the hypothesis $\|X\|_{L^2(\B_{\operatorname{max}})} \leq 1$, 
and the fact that $X \mapsto \E(X|\B)$ is an orthonormal projection we observe the estimate
\begin{equation}\label{eb}
0 \leq {\mathcal E}(\B) \leq 1.
\end{equation}
Also, if $\B \subseteq \B'$, then a simple application of Pythagoras's theorem yields
\begin{equation}\label{pythagoras}
 {\mathcal E}(\B') = {\mathcal E}(\B) + \| \E( X | \B' ) - \E( X | \B ) \|_{L^2(\B_{\operatorname{max}})}^2.
\end{equation}
In particular, finer $\sigma$-algebras have higher energy.

We shall prove the regularity lemma via an \emph{energy incrementation argument}.  We shall take some $\sigma$-algebras $\B_i, \B'_i$ and see if
they verify the required properties of the lemma.  If they do not, we will be able to replace some of these $\sigma$-algebras by finer $\sigma$-algebras
with slightly higher complexity and somewhat larger energy.  The bounds \eqref{eb}, \eqref{pythagoras} will be used to show that this energy incrementation cannot continue indefinitely, and when it does stop, we will establish the theorem.

The key step in the argument is the following.

\begin{lemma}[Lack of regularity implies energy increment]\label{energy-inc}  Suppose we have finite $\sigma$-algebras $\B'_i \subseteq \B_{i,\operatorname{max}}$ 
and events $A_i \in \B_{i,\operatorname{max}}$ for each $i \in I$ such that
$$ \left|\E\left( \bigl(X - \E( X | \bigvee_{i \in I} \B'_i )\bigr) \prod_{i \in I} 1_{A_i} \right)\right| > \frac{1}{F(M)}$$
for some $M > 0$.  Then if we set 
$$\B''_i := \B'_i \vee \{ \emptyset, A_i, \Omega \backslash A_i, \Omega \} \hbox{ for all } i \in I$$
(thus $\B''_i$ is the $\sigma$-algebra generated by $\B'_i$ and $A_i$), then we have the complexity increment
\begin{equation}\label{complexity-increment}
\complexity(\B''_i) \leq \complexity(\B'_i) + 1 \hbox{ for all } i \in I
\end{equation}
and the energy increment
\begin{equation}\label{energy-increment}
{\mathcal E}(\bigvee_{i \in I} \B''_i) \geq {\mathcal E}(\bigvee_{i \in I} \B'_i) + \frac{1}{F(M)^2}.
\end{equation}
\end{lemma}

\begin{proof}  The complexity increment is immediate from the definition of complexity.  As for the energy increment, observe that
$\prod_{i \in I} 1_{A_i}$ is measurable in $\bigvee_{i \in I} \B''_i$.  Thus we have
$$ \E\left( (X - \E( X | \bigvee_{i \in I} \B'_i )) \prod_{i \in I} 1_{A_i} \right)
= \E\left( (\E(X|\bigvee_{i \in I} \B''_i) - \E( X | \bigvee_{i \in I} \B'_i )) \prod_{i \in I} 1_{A_i}\right).$$
On the other hand, we clearly have $\E( (\prod_{i \in I} 1_{A_i})^2 ) \leq 1$.
Applying Cauchy-Schwarz, we conclude
$$ \left|\E\left( \bigl(X - \E( X | \bigvee_{i \in I} \B'_i )\bigr) \prod_{i \in I} 1_{A_i} \right)\right|^2
\leq \left\| \E(X|\bigvee_{i \in I} \B''_i) - \E( X | \bigvee_{i \in I} \B'_i ) \right\|_{L^2(\B_{\operatorname{max}})}^2.$$
By hypothesis, we thus have
$$ \left\| \E(X|\bigvee_{i \in I} \B''_i) - \E( X | \bigvee_{i \in I} \B'_i ) \right\|_{L^2(\B_{\operatorname{max}})}^2 \geq \frac{1}{F(M)^2}.$$
The claim now follows from \eqref{pythagoras}.
\end{proof}

We can now quickly prove Theorem \ref{measure-lemma}.  We shall run the following double-loop
algorithm to generate $\B_i$, $\B'_i$, and $M$.

\begin{itemize}
\item Step 0: Initialize $\B_i = \B'_i = \{ \emptyset,\Omega\}$ to be the trivial $\sigma$-algebra for each $i \in I$.
\item Step 1: Set $M$ to be the quantity
$$ M := \max\left( m, \max_{i \in I} \complexity(\B_i) \right).$$
Thus, for instance, the initial value of $M$ will be $m$.
\item Step 2: If \eqref{xbfine} holds, then we halt the algorithm.  Otherwise, we can apply Lemma \ref{energy-inc} to
locate $\sigma$-algebras $\B'_i \subseteq \B''_i \subseteq \B_{i,\operatorname{max}}$ for $i \in I$ obeying \eqref{complexity-increment} and
\eqref{energy-increment}.
\item Step 3: If we have
$$ {\mathcal E}(\bigvee_{i \in I} \B''_i) \leq {\mathcal E}(\bigvee_{i \in I} \B_i) + \eps^2$$
then we set $\B'_i$ equal to $\B''_i$ for each $i \in I$, and return to Step 2.  Otherwise, we set
$\B_i$ and $\B'_i$ \emph{both} equal to $\B''_i$ for each $i \in I$, and return to Step 1.
\end{itemize}

The following observations about the above algorithm are easily verified by induction:

\begin{itemize}
\item At every stage of the algorithm, we have $\B_i \subseteq \B'_i \subseteq \B_{i,\operatorname{max}}$ for all $i \in I$.
\item At every stage of the algorithm, we have
$$ {\mathcal E}(\bigvee_{i \in I} \B'_i) \leq {\mathcal E}(\bigvee_{i \in I} \B_i) + \eps^2$$
and hence by \eqref{pythagoras} we have \eqref{coarsefine}.
\item At every stage of the algorithm we have $m \leq M$ and $\complexity(\B_i) \leq M$ for all $i \in I$.
\end{itemize}

Thus, if the algorithm does halt (so that \eqref{xbfine} holds), then we will have achieved every objective of
Theorem \ref{measure-lemma}, except possibly for the upper bound $M = O_{F,\eps,E}(1)$ on $M$.  Hence the only remaining task is
to show that the algorithm does indeed halt in finite time with the required bound on $M$.

Let us first analyze the inner loop of the algorithm, which loops between Step 2 and Step 3.  At the start of this inner loop (i.e. when one enters Step 2 from Step 1),
the $\B'_i$ are equal to $\B_i$.  At each execution of this inner loop,
the energy ${\mathcal E}(\bigvee_{i \in I} \B'_i)$ increases by at least $\frac{1}{F(M)^2}$, thanks to \eqref{energy-increment}, while
the complexities $\complexity(\B'_i)$ increase by at most 1, thanks to \eqref{complexity-increment}.  On the other hand,
if the energy ${\mathcal E}(\bigvee_{i \in I} \B'_i)$ ever increases by more than $\eps^2$, then we will end the inner loop and instead
trigger the outer loop (returning from Step 3 to Step 1).  Thus for any fixed iteration of the outer loop, the inner loop can run for
at most $F(M)^2/\eps^2 + 1$ iterations, and the complexity of the $\sigma$-algebras $\B'_i$ increase by at most $F(M)^2/\eps^2 + 1$ when doing so.
In particular, the inner loop always terminates in finite time.

Now we can analyze the outer loop.  At the beginning of this loop, the $\B_i$ are equal to the trivial algebra, and $M$ is equal to $m$.
After each iteration of this outer loop, each $\B_i$ is replaced by a $\sigma$-algebra $\B''_i$ whose complexity is at most $F(M)^2/\eps^2 + 1$
more than the complexity of $\B_i$.  In particular, the complexity of the new value of $\B_i$ is at most $M + F(M)^2/\eps^2 + 1$, which causes
the new value of $M$ to be bounded by $M + F(M)^2/\eps^2 + 1$.  Also,
the energy ${\mathcal E}(\bigvee_{i \in I} \B_i)$ of $\B_i$ will increase by at least $\eps^2$.  From \eqref{eb} we thus see that the outer loop
can execute at most $\lfloor 1/\eps^2 \rfloor$.  Thus the algorithm terminates in finite time, and the final value of $M$ is bounded by the quantity
obtained by applying $\lfloor 1/\eps^2 \rfloor$ iterations of the map $M \mapsto M + F(M)^2/\eps^2 + 1$ to $m$, so in particular $M = O_{F,\eps,m}(1)$.
This completes the proof of Theorem \ref{measure-lemma}.
\endprf

\begin{remark} The doubly-iterated nature of the argument, combined with the desire for the growth function $F$ to be exponential for the application to
Theorem \ref{srl}, causes the final bounds on $M$ (and hence on $J$) to be tower-exponential in $1/\eps^C$ for some absolute constant $C$.
As discussed in \cite{gowers-sz}, this tower exponential bound cannot be significantly improved.  However, by lowering $F$ to linear or polynomial
growth one can obtain a somewhat weaker regularity lemma, but with better bounds; see \cite{komlos} for some further discussion on how one
can adjust the strength of the regularity lemma to suit one's application.  In the converse direction, we will need to increase $F$ further,
to tower-exponential or even faster, when we iterate this lemma to obtain hypergraph regularity lemmas\footnote{Basically, to obtain a satisfactory regularity control on hypergraphs, say $3$-uniform hypergraphs, one has to first apply a result such as Theorem \ref{measure-lemma} with some growth function $F^{fast}$ to approximate some $3$-uniform object by a collection of $2$-uniform $\sigma$-algebras (i.e. partitions of complete graphs into incomplete graphs).  One then applies Theorem \ref{measure-lemma} again with another growth function $F$ to approximate the atoms of those $2$-uniform $\sigma$-algebras by some $1$-uniform objects (vertex partitions). In order for the error terms to be manageable, it turns out that $F^{fast}$ has to grow much faster than $F$, in fact it must essentially be an iterated version of $F$.  See \cite{tao:hyper} for further discussion.}.  
The flexibility afforded by this additional
parameter $F$, which is not present in the usual formulation of the regularity lemma, may hopefully be useful for other applications also.
\end{remark}

\section{An entropy variant of the regularity lemma}

One can also give a variant of the above arguments, in which the $L^2$ norm is replaced by the Shannon entropy.  In particular, the energy
incrementation argument is replaced by an entropy incrementation argument, which gives the lemma a much more information-theoretic flavour than
before.  As always we fix an ambient probability space $(\Omega, \B_{\operatorname{max}}, \P)$.

\begin{definition}[Entropy]  If $\B \subset \B_{\operatorname{max}}$ is a finite $\sigma$-algebra, we define the \emph{Shannon entropy} $\H(\B)$ to be the quantity
$$ \H(\B) := \sum_A \P(A) \log_2 \frac{1}{\P(A)}$$
where $A$ ranges over all the atoms of $\B$ and we adopt the convention $0 \log \frac{1}{0} = 0$.  If $X$ is a random variable taking only finitely many values, we define $\H(X) := \H(\B_X)$, where $\B_X$ is the $\sigma$-algebra generated by $X$.  In other words
$$ \H(X) := \sum_x \P(X=x) \log_2 \frac{1}{\P(X=x)}.$$
\end{definition}

It is easy to verify that if $X$ is a Boolean variable (only taking the values $0$ and $1$), then $\H(X)$ can be at most 1.  More generally,
we have the inequality
$$ \H(\B) \leq \complexity(\B)$$
for any finite $\sigma$-algebra $\B$.  The quantity $\H(X)$ measures, roughly speaking, how much information one could learn from $X$.  It can be viewed as a more refined version of the complexity, which is less sensitive to exceptional events of small probability than the complexity is.

In the probabilistic formulation of the regularity lemma, conditional expectation played a prominent role.  In the entropy formulation,
the analogous concept is \emph{conditional entropy}.

\begin{definition}[Conditional entropy]  If $X, Y$ are random variables taking finitely many values, we define the conditional entropy
$\H(X|Y)$ by the formula
\begin{align*}
 \H(X|Y) &:= \sum_y \P(Y=y) \H(X|Y=y) \\
 &= \sum_y \P(Y=y) \sum_x \P(X=x|Y=y) \log_2 \frac{1}{\P(X=x|Y=y)}.
 \end{align*}
\end{definition}

An equivalent definition is given by the Bayes identity
$$ \H(X|Y) = \H(X,Y) - \H(Y).$$
The quantity $\H(X|Y)$ measures, roughly speaking, how much new information one could still learn from $X$ if one already knew the value of $Y$
(thus for instance $\H(X|X)$ is always zero).

Another key quantity we need is the \emph{conditional mutual information} $\I(X:Y|Z)$ of three random variables $X,Y,Z$ taking finitely many
values, defined by
$$ \I(X:Y|Z) := \H(X|Z) - \H(X|Y,Z) = \H(Y|Z) - \H(Y|X,Z);$$
informally, it measures how much knowing $Y$ would tell one about $X$, or vice versa, assuming that $Z$ is already known.
A handy (and intuitive) fact is that the conditional mutual information is always non-negative; this is equivalent to
the \emph{submodularity inequality} 
$$ \H(X,Y,Z) + \H(Z) \leq \H(X,Z) + \H(Y,Z)$$
for entropy, and can be proven via Jensen's inequality.  A more quantitative assertion of this fact is given in
Lemma \ref{entropy-expectation} below.

If $X$ and $Y$ are random variables, we write $X \mapsto Y$, and say that $Y$ is \emph{determined} by $X$, if $\B_Y \subseteq \B_X$.
If $X$ and $Y$ take only finite values, then $X \mapsto Y$ is equivalent to the existence of a functional relationship $Y = f(X)$ for some
deterministic function $f$, and is also equivalent (up to events of probability zero) to the conditional entropy $\H(Y|X)$ vanishing.

We now give the information-theoretic analogue of Theorem \ref{measure-lemma}.  To simplify the notation a little bit we will
restrict to the case $I= \{1,2\}$, although the generalization to more than two reference $\sigma$-algebras is not difficult.

\begin{lemma}[Information-theoretic regularity lemma]\label{mainlemma}  Let $X_1, X_2, Y$ be random variables taking finitely many values such that
$\H(Y) \leq m$ for some $m \geq 0$.  Let $F: \R^+ \to \R^+$ be an arbitrary function, and $\eps > 0$.  Then there exists random variables $Z_1, Z_2$ (the ``coarse approximation'') and $Z'_1, Z'_2$ (the ``fine approximation''), also taking finitely many values, with 
the following properties.

\begin{itemize}

\item (Determinism) We have the determinism relations 
\begin{equation}\label{determin}
X_1 \mapsto Z'_1 \mapsto Z_1; \quad X_2 \mapsto Z'_2 \mapsto Z_2.
\end{equation}

\item (Coarse approximation has bounded entropy)  We have 
\begin{equation}\label{coarse-low}
\H(Z_1, Z_2) \leq \H(Z'_1,Z'_2) = O_{F,\eps, m}(1).
\end{equation}

\item (Coarse and fine approximations are close) We have
\begin{equation}\label{coarse-fine}
\I( Y : Z'_1, Z'_2 | Z_1, Z_2 ) \leq \eps.
\end{equation}

\item (Fine approximation is nearly optimal)  For any random variables $W_1, W_2$ with $X_1 \mapsto W_1$ and $X_2 \mapsto W_2$ we have
\begin{equation}\label{fine-optimal}
\I( Y : W_1,W_2 | Z'_1,Z'_2 ) \leq \frac{\H(W_1,W_2)}{F(\H(Z_1, Z_2))}.
\end{equation}

\end{itemize}

\end{lemma}

\begin{proof}
To construct $Z_1, Z_2, Z'_1, Z'_2$ we perform the following ``entropy incrementation'' algorithm, which is closely analogous
to the energy incrementation algorithm used in the proof of Theorem \ref{measure-lemma}.

\begin{itemize}
\item Step 0. Initialize $Z_1 = Z_2 = 0$ (one can of course replace $0$ by any other deterministic random variable).
\item Step 1. Let $Z'_1, Z'_2$ be random variables which minimize the quantity
\begin{equation}\label{emax}
\H(Y | Z'_1, Z'_2 ) + \frac{\H(Z'_1,Z'_2)}{F(\H(Z_1,Z_2))} 
\end{equation}
subject to the constraints $X_1 \mapsto Z'_1 \mapsto Z_1$ and $X_2 \mapsto Z'_2 \mapsto Z_2$.
(If there are several such minimizers, we select among them arbitrarily.)
\item Step 2.  If we have
$$\H(Y|Z_1,Z_2) - \H(Y|Z'_1,Z'_2) > \eps$$
then we replace $Z_1, Z_2$ with $Z'_1$, $Z'_2$ respectively, and return to Step 1.  Otherwise,
we terminate the algorithm.
\end{itemize}

We remark that because $X_1, X_2$ take only finitely many values, the number of possibilities for the random variables $Z'_1,Z'_2$ 
is finite up to equivalence. Hence a minimizer to the quantity \eqref{emax} always exists.  
Intuitively, $Z'_1,Z'_2$ is constructed to capture as much information about
$Y$ as is possible while remaining determined by $X_1,X_2$; the slight penalty term in \eqref{emax} is designed to keep some control of the 
entropy of $Z'_1,Z'_2$ (otherwise it would be as large as that of $X_1,X_2$, for which we have no bounds).  
Observe that every time we return from Step 2 to
Step 1, the quantity $\H(Y|Z_1,Z_2)$ (which measures the amount of information in $Y$ that remains to be captured by $Z_1,Z_2$) decreases by at least $\eps$.  On the other hand, from Jensen's inequality one can verify that
$$ 0 \leq \H(Y|Z_1,Z_2) \leq \H(Y) \leq m.$$
Thus the above algorithm must halt after at most $m/\eps$ iterations.  It is also clear that the random variables $Z_1, Z_2, Z'_1, Z'_2$ generated by this algorithm will obey the determinism relationships \eqref{determin} and \eqref{coarse-fine}.

Also, if $W_1, W_2$ are any random variables determined by $X_1, X_2$ respectively, then by comparing the minimizer $Z'_1, Z'_2$
against the competitor $(Z'_1,W_1)$, $(Z'_2,W_2)$ (which obeys the required constraints), we have
$$
\H(Y | Z'_1, Z'_2 ) + \frac{\H(Z'_1,Z'_2)}{F(\H(Z_1,Z_2))} 
\leq \H(Y | Z'_1, Z'_2, W_1, W_2 ) + \frac{\H(Z'_1,Z'_2, W_1, W_2)}{F(\H(Z_1,Z_2))}.$$
Since $\H(Y | Z'_1, Z'_2 ) - \H(Y | Z'_1, Z'_2, W_1, W_2 ) = \I( Y : W_1, W_2 | Z'_1, Z'_2 )$ and
$\H(Z'_1,Z'_2, W_1, W_2) \leq \H(Z'_1,Z'_2) + \H(W_1,W_2)$, we obtain \eqref{fine-optimal} as desired after some algebra.

Now we compare the entropies of $Z_1,Z_2$ and $Z'_1, Z'_2$.  
Since $Z_1, Z_2$ obeys the constraints in the minimization problem \eqref{emax}, we have
$$
\H(Y | Z'_1, Z'_2 ) + \frac{\H(Z'_1,Z'_2)}{F(\H(Z_1,Z_2))} 
\leq \H(Y | Z_1, Z_2 ) + \frac{\H(Z_1,Z_2)}{F(\H(Z_1,Z_2))}.$$
As observed earlier, the first summand on either side ranges between 0 and $m$.  Thus we have (after some rearranging)
$$ \H(Z'_1, Z'_2) \leq \H(Z_1,Z_2) + m F(\H(Z_1,Z_2)).$$
In particular, every time we return from Step 2 to Step 1, the quantity $\H(Z_1,Z_2)$ increases by at most 
$m F(\H(Z_1,Z_2))$.
From Step 0, the initial value of $\H(Z_1,Z_2)$ is 0.  Since the number of iterations is bounded by
$m/\eps$, we see that the final value of $\H(Z_1,Z_2)$ is bounded by a finite
(but extremely large) quantity $O_{m,F,\eps}(1)$ or more explicitly the value obtained after
$m/\eps$ iterations of the map $M \mapsto M +  m F(M)$ applied to $0$.
\end{proof}

To pass from an entropy formulation to an expectation formulation, we need a way to pass from control of entropy to control of expectations.
A clue to how to do this is provided by the following observation: if $Y \mapsto Y'$ and $\I(X:Y|Y') = 0$, then 
$X$ and $Y$ are independent conditionally on $Y'$.  In particular, if $X$ takes values in 
a vector space, this implies that $\E(X|Y) = \E(X|Y')$.  In other words, whenever $\I(X:Y|Y') = \H(X|Y')-\H(X|Y)$ is zero, so is $\E(X|Y') -\E(X|Y)$.  This may help motivate the following lemma, which is a perturbative version of the above observation.

\begin{lemma}[Relation between entropy and expectation]\label{entropy-expectation}
Let $X, Y, Y'$ be discrete random variables with $Y \mapsto Y'$, and with $X$ taking values in the unit interval $\{-1 \leq x \leq 1\}$.
Then we have
$$ \E\left( \bigl|\E(X|Y') - \E(X|Y)\bigr| \right) \leq 2 \I(X:Y|Y')^{1/2}.$$
\end{lemma}

More informally, this lemma asserts that approximate conditional independence in the entropy sense implies approximate conditional
independence in an expectation sense.  The bound $2 \I(X:Y|Y')^{1/2}$ is not best possible, but any bound which decays
to zero as $\I(X:Y|Y') \to 0$ will be sufficient for our purposes.

\begin{proof}  The basic idea is to exploit the observation that the function $x\log \frac{1}{x}$ is not only concave but also strictly concave
on $[0,1]$.  Let us first verify the lemma in the special case when $Y'$ is deterministic (so the hypothesis $Y \mapsto Y'$ is vacuous), thus we wish to prove
$$ \E( |\E(X) - \E(X|Y)| ) \leq 2 \I(X:Y)^{1/2}.$$
Let $1 \leq x_1,\ldots,x_n \leq -1$ be the essential range of $X$, and let $y_1,\ldots,y_m$ be the essential range of $Y$.  For any $1 \leq i \leq n$
and $1 \leq j \leq m$, define the probabilities
\begin{align*}
p_{ij} &:= \P( X = x_i | Y = y_j)\\
q_j := \P(Y = y_j)\\
\overline{p_i} &:= \sum_{j=1}^m q_j p_{ij} = \P(X = x_i)
\end{align*}
Then we observe that $0 \leq p_{ij}, q_j \leq 1$ and that $\sum_{j=1}^m q_j = 1$.  If we define $f: [0,1] \to \R$ to be
the function $f(x) := -x \log x$ (with the convention $f(0) := 0$), we thus have 
\begin{align*}
\I(X:Y) &= \H(X) - \H(X|Y) \\
&= \sum_{i=1}^n (f(\overline{p_i}) - \sum_{j=1}^m q_j f(p_{ij})).
\end{align*}
Now observe that $f$ is concave, indeed we have $f''(x) = -1/x$ for all $x \in (0,1]$.  Thus by Taylor's theorem with remainder,  
$$ f(p_{ij}) \leq f(\overline{p_i}) + f'(\overline{p_i}) (p_{ij} - \overline{p_i}) - \frac{1}{2} (p_{ij} - \overline{p_i})^2 / p^*_{ij}$$
where $p^*_{ij}$ is a quantity between $p_{ij}$ and $\overline{p_i}$.
Inserting this into the preceding estimate and noting that $\sum_{j=1}^m q_j(p_{ij} - \overline{p_i})  = 0$, we conclude
that
$$ \sum_{j=1}^m q_j \sum_{i=1}^n (p_{ij} - \overline{p_i})^2 / p^*_{ij} \leq 2 \I(X:Y).$$
Now we compute using the boundedness of $x_i$ and Cauchy-Schwarz, as well as the crude estimate $p^*_{ij} \leq \overline{p_i} + p_{ij}$,
\begin{align*}
\E( |\E(X) - \E(X|Y)| ) &= \sum_{j=1}^m q_j |\E(X) - \E(X|Y = y_j)| \\
&=\sum_{j=1}^m q_j |\sum_{i=1}^n x_i (\overline{p_i} - p_{ij})| \\
&\leq \sum_{j=1}^m q_j \sum_{i=1}^n  |\overline{p_i} - p_{ij}| \\
&\leq (\sum_{j=1}^m q_j \sum_{i=1}^n |\overline{p_i} - p_{ij}|^2/p^*_{ij})^{1/2}
(\sum_{j=1}^m q_j \sum_{i=1}^n p^*_{ij})^{1/2} \\
&\leq [2 \I(X:Y) \sum_{j=1}^m q_j \sum_{i=1}^n \overline{p_i} + p_{ij}]^{1/2}\\
&= 2 \I(X:Y)^{1/2}.
\end{align*}
Now we consider the general case when $Y'$ is not deterministic.  In that case we write
$$ \E( |\E(X|Y') - \E(X|Y)| ) = \sum_{y'} \P(Y'=y')
\E( |\E(X|Y'=y') - \E(X|Y;Y'=y')| ).$$
(Here we have taken advantage of the hypothesis $Y \mapsto Y'$.)  Applying the preceding computation, we conclude
$$ \E( |\E(X|Y') - \E(X|Y)| ) \leq \sum_{y'} \P(Y'=y')
2 \I(X:Y|Y'=y')^{1/2}.$$
Applying Cauchy-Schwarz again we conclude
\begin{align*}
\E( |\E(X|Y') - \E(X|Y)| ) &\leq 2\sqrt{\sum_{y'} \P(Y'=y')
\I(X:Y|Y'=y') } \\
&= 2 \I(X:Y|Y')^{1/2}
\end{align*}
as desired.
\end{proof}

By combining this with Lemma \ref{mainlemma} it is possible to give a statement closely resembling Theorem \ref{measure-lemma}, and which is
also sufficient to imply Theorem \ref{srl}.  We omit the details.

\end{document}